\documentclass[12pt]{article}

\def\b#1{{\bf #1}}
\def\c#1{{\cal #1}}
\def\1{{\bf 1}}
\newcommand{\botimes}{\mbox{\boldmath \boldmath$\otimes$}}
\newcommand{\tl}{\,\triangleleft}
\def\cross{{\triangleright\!\!\!<}}
\def\cocross{{>\!\!\!\triangleleft\,}}
\def \uqg{\mbox{$U_q{\/\mbox{\bf g}}$ }}
\def\R{{\cal R}\,}

\newtheorem{prop}{Proposition}
\newtheorem{theorem}{Theorem}

\newcommand{\be}{\begin{equation}}
\newcommand{\ee}{\end{equation}}
\newcommand{\bea}{\begin{eqnarray}}
\newcommand{\eea}{\end{eqnarray}}

\begin{document}

\title{Decoupling of Translations from Homogeneous Transformations
in Inhomogeneous Quantum Groups\footnote{Talk given at the 2nd International Symposium on
``Quantum Theory and Symmetries'', July  2001, Krakow, Poland.
Preprint 00-45 Dip. Matematica e Applicazioni, Universit\`a di Napoli;
Preprint DSF/36-2001.}}

\author{        Gaetano Fiore, \\\\
         \and
        Dip. di Matematica e Applicazioni, Fac.  di Ingegneria\\ 
        Universit\`a di Napoli, V. Claudio 21, 80125 Napoli
        \and
        I.N.F.N., Sezione di Napoli,\\
        Complesso MSA, V. Cintia, 80126 Napoli
        }
\date{}

\maketitle

\abstract{We briefly report on our result that, 
if there exists a realization of a Hopf
algebra $H$ in a $H$-module algebra $\c{A}$, then 
their cross-product is equal
to the product of $\c{A}$ itself
with a subalgebra isomorphic to $H$ and {\it commuting} with 
$\c{A}$. We illustrate its application to the Euclidean
quantum groups in $N\ge 3$ 
dimensions.}

\section{Introduction and main results}

Let $H$ be a Hopf algebra and  $\c{A}$ a unital $H$-module algebra.
Here for the sake of brevity
we stick to consider the case of a right $H$-module algebra.
Then by definition there exists a bilinear map
$\tl: (a,g)\in\c{A} \times H\to a\tl g\in\c{A}$, 
called a right action, such that
\bea
&&a\tl(gg') = (a\tl g) \tl   g'                  \label{modalg1}\\[4pt]
&&(aa')\tl g = (a\tl g_{(1)})\, (a'\tl g_{(2)}).       \label{modalg2}
\eea
We use a Sweedler-type notation
with suppressed summation sign for the coproduct $\Delta(g)$ 
of $g$, namely the short-hand notation 
$\Delta(g)=g_{(1)}\botimes g_{(2)}$. 
The cross-product algebra $\c{A} \cocross H$ is $H\otimes\c{A}$ 
as a vector space, whereas
the product of two elements is given by
\be
(g\otimes a)(g'\otimes a' )=gg'_{(1)} \otimes (a\tl g'_{(2)})a'
                                            \label{product}
\ee
for any $a,a'\in\c{A}$, $g,g'\in H$. To simplify the notation, 
in the sequel we denote 
$g\otimes a$ by $ga$ and omit either unit $\1_{\c{A}},\1_H$ 
whenever multiplied
by non-unit elements; consequently, for $g=\1_H$, $a'=\1_{\c{A}}$
(\ref{product}) takes the form
\be
a g'=g'_{(1)}\, (a\tl g'_{(2)}).                \label{crossprod}
\ee
$\c{A}\cocross H$ itself is a $H$-module algebra under the right
action  $\tl$ of $ H$ if we extend the latter on the elements of $H$
as the adjoint action,
\be
h\tl g=Sg_{(1)} h g_{(2)},          \qquad\qquad g,h\in  H. \label{adjor}
\ee
One very important situation in which one encounters cross products
is in inhomogeneous (quantum) groups, that are Hopf algebras
whose algebras are cross-products $\c{A}\cocross H$, where $H$
is the Hopf subalgebra of ``homogeneous transformations'' and $\c{A}$ 
the (braided) Hopf subalgebra of ``translations''.

In Ref. \cite{Fio01} we have found a sufficient condition for 
$\c{A}\cocross H$ to be equal to a product $\c{A}\, H'$,
where $H'\subset H \cross \c{A}$ is a subalgebra 
isomorphic to $H$ and {\it commuting} 
with $\c{A}$: it is sufficient 
to assume that there exists an algebra map
\be
\tilde\varphi: \c{A}\cocross H \rightarrow \c{A}   \label{Homr}
\ee
acting as the identity on the subalgebra
$\c{A}$, namely for any $a\in \c{A}$
\be
\tilde\varphi(a)=a .                          \label{ident0r}
\ee
Here we briefly recall this result, as well as a similar one
requiring a weaker assumption, namely that $H$ admits a Gauss
decomposition into two Hopf subalgebras $H^+,H^-$ for each
of which analogous maps $\tilde\varphi^+,\tilde\varphi^-$
exist. In the last section we illustrate the application
of these general results to a class of inhomogeneous quantum
groups, the Euclidean quantum groups in $N\ge 3$ dimensions,
for which $\c{A}$ is (``the algebra of functions on'')
the quantum Euclidean $N$-dimensional space; the
corresponding maps  $\tilde\varphi,\tilde\varphi^+,\tilde\varphi^-$
have been determined in \cite{CerFioMad00} and further
studied in \cite{FioSteWes00}. 
[Note that no such maps can exist for the (undeformed) Euclidean algebra, 
for which $\c{A}\equiv$the algebra of functions on
$\b{R}^N$, which is abelian].
Other applications
include \uqg-covariant Heisenberg (or Clifford) algebras
and the $q$-deformed 2-dimensional fuzzy sphere $S_{q,M}^2$ 
\cite{fuzzyq}.

\medskip

Let $\tilde\c{C}$ be the commutant of $\c{A}$ within $\c{A}\cocross H$, 
i.e. the subalgebra 
\be
\tilde\c{C}:=\{c\in \c{A}\cocross H \: |\:
[c,a]=0  \: \:\: \: \forall a\in\c{A}\}.                \label{commu}
\ee
Clearly $\tilde\c{C}$ contains the center $\c{Z}(\c{A})$ of $\c{A}$.
Let $\tilde\zeta: H\rightarrow \c{A}\cocross H$ be the map defined by
\be
\tilde\zeta(g):=g_{(1)}) \tilde\varphi(S g_{(2)}).    \label{deftildezeta}
\ee
Note that if we apply $\tilde\varphi$ to $\tilde\zeta$ and 
recall that $\tilde\varphi$ is both
a homomorphism and idempotent we find
$\tilde\varphi\circ\tilde\zeta=\varepsilon$.
Now, $\varepsilon(g)$ is a complex number times $\1_{\c{A}}$
and therefore trivially commutes with $\c{A}$. Since by definition
the commutation relations between $\c{A}$ and either 
$H$ or $\tilde\varphi(H)$ are the same, we expect that
also  $\tilde\zeta(g)$ commutes with $\c{A}$. This is confirmed by

\begin{theorem}\cite{Fio01} Under the above assumptions
the map $\tilde\zeta$ is an injective algebra homomorphism
$\tilde\zeta: H\rightarrow\tilde\c{C}$; moreover 
$\tilde\c{C}=\tilde\zeta(H)\,\c{Z}(\c{A})$ and
$\c{A}\cocross H=\tilde\zeta(H)\,\c{A}$.
If, in particular  $\c{Z}(\c{A})=\b{C}$, then $\tilde\c{C}=\zeta(H)$
and $\tilde\zeta: H\leftrightarrow\tilde\c{C}$ is an algebra isomorphism.
Moreover, the center of the cross-product $\c{A}\cocross H$ is given by
$\c{Z}(\c{A}\cocross H)=\c{Z}(\c{A})\,\tilde\zeta\left(\c{Z}(H)\right)$.
Finally, if $H_c,\c{A}_c$ are maximal abelian subalgebras of $H$
and $\c{A}$ respectively, then $\c{A}_c\,\tilde\zeta(H_c)$
is a maximal abelian subalgebra of $\c{A}\cocross H$.
\label{theorem1r}
\end{theorem}

In other words, the subalgebra $H'$ looked for
will be obtained by setting $H':=\tilde\zeta(H)$.
For these reasons we shall call $\tilde\zeta$, as well as the other
maps $\tilde\zeta^{\pm}$ which we shall introduce below,
{\it decoupling maps}.

\medskip
To introduce the latter we just need a weaker assumption, namely
that, instead of a map $\tilde\varphi$, we just have 
at our disposal two homomorphisms $\tilde\varphi^+,\tilde\varphi^-$ 
\be
\tilde\varphi^{\pm}: H^{\pm} \cross \c{A} \rightarrow \c{A}      \label{Hom+-}
\ee
fulfilling (\ref{ident0r}), 
where $H^+,H^-$ denote two Hopf subalgebras of $H$ such that
Gauss  decompositions $H=H^+H^-=H^-H^+$ hold.
(The typical case is when $H=\uqg$ and $H^+,H^-$ denote its
positive and negative Borel subalgebras.)
Then Theorem \ref{theorem1r} will apply separately to 
$\c{A}\cocross H^+$ and $\c{A}\cocross H^-$, if we define
corresponding maps $\tilde\zeta^{\pm}: H^{\pm} \rightarrow \c{A}$ by
\be
\tilde\zeta^{\pm}(g):=g_{(1)}\tilde\varphi^{\pm}(S  g_{(2)}), \label{defzeta+-r}
\ee
where $g\in H^{\pm}$ respectively. What can we say about
the whole $H\cross \c{A}$? 

\begin{theorem}\cite{Fio01} Under the above assumptions formulae
(\ref{defzeta+-r}) define injective algebra homomorphisms
$\tilde\zeta^{\pm}: H^{\pm}\rightarrow\tilde\c{C}$. Moreover,
\be 
\tilde\c{C}=\tilde\zeta^+(H^+)\,\tilde\zeta^-(H^-)\,\c{Z}(\c{A})
=\tilde\zeta^-(H^-)\,\tilde\zeta^+(H^+)\,  \c{Z}(\c{A})       \label{casim} 
\ee 
and
\be
\c{A} \cocross H=\tilde\zeta^+(H^+)\,\tilde\zeta^-(H^-)\,\c{A}
=\tilde\zeta^-(H^-)\,\tilde\zeta^+(H^+)\,\c{A}.         \label{decom}
\ee
In particular, if  $\c{Z}(\c{A})=\b{C}$, 
then the commutant is equal to
$\tilde\c{C}=\tilde\zeta^+(H^+)\tilde\zeta^-(H^-)=
\tilde\zeta^-(H^-)\tilde\zeta^+(H^+)$.
Furthermore, any element $c$ of the center $\c{Z}(\c{A}\cocross H)$
of the cross-product $\c{A}\cocross H$
can be expressed in the form 
\be
c=\tilde\zeta^+\left(c^{(1)}\right)\tilde\zeta^-\left(c^{(2)}\right)c^{(3)},
\ee
where $c^{(1)}\otimes c^{(2)}\otimes c^{(3)}
\in H^+\otimes H^- \otimes\c{Z}(\c{A})$ and
$c^{(1)}c^{(2)}\otimes c^{(3)}
\in \c{Z}(H) \otimes\c{Z}(\c{A})$; vice versa any
such object $c$ is an element of $\c{Z}(\c{A}\cocross H)$. Finally,
if $H_c\subset H^+ \cap H^-$ and 
$\c{A}_c$ are  maximal abelian subalgebras of $H$
and $\c{A}$ respectively, then $\c{A}_c\,\tilde\zeta^+(H_c)$
(as well as $\c{A}_c\,\tilde\zeta^-(H_c)$)
is a maximal abelian subalgebra of $\c{A}\cocross H$.
\label{theorem+-r}
\end{theorem}

As a consequence of this theorem, for any $g^+\in H^+$, $g^-\in H^-$
there exists a sum $c^{(1)}\otimes c^{(2)}\otimes c^{(3)}\in 
\c{Z}(\c{A})\otimes H^-\otimes H^+$
(depending on $g^+,g^-$) such that
\be
\tilde\zeta^+(g^+)\tilde\zeta^-(g^-)=c^{(1)}\tilde
\zeta^-(c^{(2)})\tilde\zeta^+(c^{(3)}).            \label{com+-}
\ee
These will be the commutation relations between elements of
$\tilde\zeta^+(H^+)$ and $\tilde\zeta^-(H^-)$. Their form will depend on
the specific algebras considered.

\medskip

Assume that $H$ is a Hopf $*$-algebra and $\c{A}$ a $H$-module
$*$-algebra. Then, as known, 
these two $*$-structures can be glued in a unique one to make 
$\c{A}\cocross H$ a $*$-algebra itself. What can we say about
the behaviour of the decoupling maps under the latter 
$*$-structure?

\begin{prop}\cite{Fio01}
If $\tilde\varphi: \c{A}\cocross H \rightarrow \c{A}$
is a $*$-homomorphism, then also the map
$\tilde\zeta: H \rightarrow \tilde\c{C}$ is.
If $\tilde\varphi^{\pm}$ are $*$-homomorphism, then also the map
$\tilde\zeta^{\pm}: H^{\pm} \rightarrow \tilde\c{C}$ are. 
If $\tilde\varphi^{\pm}$
fulfill the condition
$\tilde\varphi^{\pm}(\alpha^*)=[\tilde\varphi^{\mp}(\alpha)]^*$
for any $\alpha\in H^{\mp} \cross \c{A}$,
then $\tilde\zeta^{\pm}$ fulfill
\be
\tilde\zeta^{\pm}(g^*)=[\tilde\zeta^{\mp}(g)]^*,\qquad\qquad g\in
H^{\mp}.                              
\ee
\label{*prop}
\end{prop}

\section{Application to the Euclidean quantum 
group $\b{R}_q^N\cocross U_qso(N)$}
\label{applications}

As an algebra $\c{A}$ we shall consider a slight extension of
the quantum Euclidean space $\b{R}_q^N$
~\cite{FadResTak89} (the $U_qso(N)$-covariant quantum space),
i.e.  of the unital associative algebra generated by $p^i$ fulfilling the
relations 
\be 
\c{P}_a{}^{ij}_{hk}p^hp^k=0, \label{xxrel} 
\ee 
where $\c{P}_a$
denotes the $q$-deformed antisymmetric projector (\ref{projectorR'}),
and the indices take the values
$i=-n,\ldots,-1,0,1,\ldots n$ for $N$ odd,
and $i=-n,\ldots,-1, 1,\ldots n$ for $N$ even; here
$n:=\left[\frac N 2\right]$ is the rank of $so(N)$, .
The multiplet $(p^i)$ carries the fundamental $N$-dim
(or vector) representation $\rho$
of $U_qso(N)$: for any $g\in U_qso(N)$
\be 
p^i\tl g=\rho^i_j(g)p^j.  \label{fund1} 
\ee
As a set of generators of $H\equiv U_qso(N)$
it is convenient to introduce the FRT generators \cite{FadResTak89}
$\c{L}^+{}_j^i,\c{L}^-{}_j^i$,  together with
the square roots of the
diagonal elements $\c{L}^+{}_j^j,\c{L}^-{}_j^j$.
The FRT generators and the braid matrix $\hat R$
are related to the so-called universal $R$-matrix 
$\R\in H\otimes H$ by
\be
\c{L}^+{}_l^a:=\R^{(1)}\rho_l^a(\R^{(2)})\qquad\qquad
\c{L}^-{}_l^a:=\rho_l^a(\R^{-1}{}^{(1)})\R^{-1}{}^{(2)} \label{frt}
\ee
and $\hat R^{ij}_{hk}=\rho^j_h(\c{L}^+{}_k^i)=(\rho^j_h\otimes\rho^i_k)(\R)$.
Since in our conventions $\R\in H^+\botimes H^-$ 
($H^+,H^-$ denote the positive, negative Borel subalgebras) we see that 
$\c{L}^+{}_l^a\in H^+$ and $\c{L}^-{}_l^a\in H^-$. 

To define $\tilde\varphi$ or 
$\tilde\varphi^{\pm}$ one \cite{CerFioMad00} slightly enlarges 
$\b{R}_q^N$ as follows. One 
introduces some new generators $\sqrt{P_a}$, with
$1\le a\le \frac N 2$, together with their  
inverses $(\sqrt{P_a})^{-1}$, requiring that 
\be
P_a^2=\sum\limits_{h=-a}^a p^hp_h=\sum\limits_{h,k=-a}^a g_{hk}p^hp^k.
\ee
In the previous equation $g_{hk}$ denotes the `metric matrix' of 
$SO_q(N)$:
\be
g_{ij}=g^{ij}=q^{-\rho_i} \delta_{i,-j}.          \label{defgij}
\ee
It is a $SO_q(N)$-isotropic tensor and is a deformation of the 
ordinary Euclidean metric.
We have introduced the notation \cite{FadResTak89}
$(\rho_i)=(n-\frac{1}{2},\ldots,\frac{1}{2},0,-\frac{1}{2},
\ldots,\frac{1}{2}-n)$
for $N$ odd, $(n-1,\ldots,0,0,\ldots,1-n)$ for $N$ even. 
In the sequel we shall call $P_n^2$ also $P^2$.
Moreover for odd $N$ we add also $\sqrt{p^0}$ and its inverse
as new generators. The commutation relations 
involving these new generators can be fixed consistently, in 
particular one finds that
\be
\sqrt{P_a}p^i=p^i\sqrt{P_a}\times\cases{1 \mbox{ if } |i|\le a \cr
\sqrt{q}\mbox{ if }i> a \cr 
1/\sqrt{q}\mbox{ if } i< a \cr}
\ee
The center of $\b{R}_q^N$ is generated by $\sqrt{P}$ and, only in the case
of even $N$, by $\sqrt{p^1/p^{-1}}$ and its inverse
$\sqrt{p^{-1}/p^1}$.
In the case of even $N$ one needs to include also
the FRT generator  $\c{L}^-{}^1_1=\c{L}^+{}^{-1}_{-1}$ 
and its inverse $\c{L}^+{}^1_1=\c{L}^-{}^{-1}_{-1}$
[which are generators of $U_qso(N)$ belonging to the natural Cartan 
subalgebra] among the generators of $\c{A}$. 
They satisfy nontrivial commutation relations
with the generators of $\c{A}$, and the standard FRT relations with the rest
of $U_qso(N)$. As a consequence, 
$\sqrt{p^{\pm 1}/p^{\mp1}}$ are eliminated from the center of
$\c{A}$ (in fact they do not $q$-commute with $\c{L}^{\pm}{}^1_1$).

The homomorphisms \cite{CerFioMad00} 
$\tilde\varphi^{\pm}:\c{A} \cocross  U_q^{\pm}so(N) \rightarrow \c{A}$
take the simplest and most compact expression on the FRT generators of
$U_q^{\pm}so(N)$. Let us introduce the short-hand notation
$[A,B]_x=AB-xBA$.     
The images of $\tilde\varphi^-$ on the FRT generators
read
\bea
&&\tilde\varphi^-(\c{L}^-{}^i_j)=g^{ih}[\mu_h,p^k]_qg_{kj},  
\label{imagel-}\\[4pt]
&&\tilde\varphi^+(\c{L}^+{}^i_j)=g^{ih}[\bar\mu_h,p^k]_{q^{-1}}g_{kj},
                                                     \label{imagel+}
\eea
where
\be
\begin{array}{lll}
\mu_0=\gamma_0 (p^0)^{-1} &\qquad\bar\mu_0=\bar\gamma_0 (p^0)^{-1}
&\quad\mbox{for $N$ odd,} \\[4pt]
\mu_{\pm 1}=\gamma_{\pm 1} (p^{\pm 1})^{-1} \c{L}^{\pm}{}^1_1
 &\qquad \bar\mu_{\pm 1} = 
\bar\gamma_{\pm 1} (p^{\pm 1})^{-1} \c{L}^{\mp}{}^1_1
&\quad\mbox{for $N$ even,} \\[4pt]
\mu_a=\gamma_a P_{|a|}^{-1}P_{|a|-1}^{-1} p^{-a} &\qquad 
\bar\mu_a = \bar\gamma_a P_{|a|}^{-1}P_{|a|-1}^{-1} p^{-a}
&\quad\mbox{otherwise,} 
\end{array}                                             \label{defmu}
\ee
and $\gamma_a, \bar\gamma_a \in \b{C}$ are normalization constants fulfilling
the conditions
\be
\begin{array}{lll}
\gamma_0 = -q^{-\frac{1}{2}} h^{-1} &\qquad 
\bar \gamma_0 = q^{\frac{1}{2}} h^{-1} &\quad\mbox{for $N$ odd,} \\[4pt]
\gamma_{\pm 1}=-k^{-1} &\qquad \bar \gamma_{\pm 1}=k^{-1} 
&\quad\mbox{for $N$ even,}\\[4pt]
\gamma_1 \gamma_{-1}=-q^{-1} h^{-2}  &\qquad \bar \gamma_1 \bar \gamma_{-1} 
= -q h^{-2}&\quad\mbox{for $N$ odd,} \\[4pt]
\gamma_a \gamma_{-a} =
-q^{-1} k^{-2} \omega_a \omega_{a-1} &\qquad \bar \gamma_a \bar \gamma_{-a} 
= -q k^{-2} \omega_a \omega_{a-1} &\quad\mbox{for $a>1$}. \nonumber
\end{array}                                               \label{gamma}
\ee
Here $h:=q^{\frac{1}{2}}-q^{-\frac{1}{2}}$,
$k:=q\!-\!q^{-1}$, $\omega_a:=(q^{\rho_a}+q^{-\rho_a})$.
One can choose the $\gamma_a$'s, $\bar\gamma_a$'s: 1) only for odd $N$, in such
a way that $\tilde\varphi^-,\tilde\varphi^+$ can be ``glued'' into a
unique homomorphism $\varphi$ of the type (\ref{Homr}), (\ref{ident0r});
2) for $|q|=1$ or $q\in\b{R}^+$ respectively, 
in such a way that either condition necessary
for the application of Proposition
\ref{*prop} is fulfilled.
{}From definition (\ref{deftildezeta}), 
using (\ref{coprodL}), (\ref{imagel-}), (\ref{imagel+}), 
we find
\be
\begin{array}{l}
\tilde\zeta^-(\c{L}^-{}^i_j)=\c{L}^-{}^i_h\tilde\varphi^-(S\c{L}^-{}^h_j)=
\c{L}^-{}^i_h[\mu_h,p^i]_q, \\[4pt]
\tilde\zeta^+(\c{L}^+{}^i_j)=\c{L}^+{}^i_h\tilde\varphi^+(S\c{L}^+{}^h_j)=
\c{L}^+{}^i_h[\bar\mu_h,p^i]_{q^{-1}}. 
\end{array}          \label{z+-}
\ee
In Ref. \cite{Fio01} 
we have determined the commutation 
relations (\ref{com+-}) between them.

\section*{Appendix}

The braid matrix of $U_qso(N)$ is given by
\bea
\hat R&=&q \sum_{i \neq 0} e^i_i \otimes e^i_i +
\sum_{\stackrel{\scriptstyle i \neq j,-j} 
{\mbox{ or } i=j=0}} e^j_i \otimes e^i_j+ q^{-1} 
\sum_{i \neq 0} e^{-i}_i
\otimes e^i_{-i} 
\label{defRsoN} \\[4pt]
&&+k (\sum_{i<j} e^i_i \otimes e^j_j- 
\sum_{i<j} q^{-\rho_i+\rho_j} 
e^{-j}_i \otimes e^j_{-i}). \nonumber
\eea
It admits the orthogonal projector
decomposition
\be
\hat R = q\c{P}_s - q^{-1}\c{P}_a + q^{1-N}\c{P}_t;      
\label{projectorR'}
\ee
$\c{P}_a,\c{P}_t,\c{P}_s$ are $q$-deformed antisymmetric, trace,
trace-free symmetric projectors. 
The coproduct and antipode of the FRT generators are given by
\be
\Delta(\c{L}^{\pm}{}^i_j)=\c{L}^{\pm}{}^i_h\otimes\c{L}^{\pm}{}^h_j,
\qquad \qquad S\c{L}^{\mp}{}^j_i=g_{ih}\c{L}^{\mp}{}^h_k g^{kj}. 
\label{coprodL}
\ee


\begin{thebibliography}{99}

\bibitem{CerFioMad00}
B. L. Cerchiai, G. Fiore, J. Madore,
`` Geometrical Tools for Quantum Euclidean Spaces'',
{\em Commun. Math. Phys.} {\bf 217} (2001), 521-554.

\bibitem{FadResTak89}
L.D.~Faddeev, N.Y.~Reshetikhin, L.~Takhtadjan, 
{\em Alge. i Analy.} {\bf 1} (1989), 178, translated from the
Russian in {\em Leningrad Math. J.} {\bf 1} (1990), 193.

\bibitem{Fio01}
G. Fiore, {\em On the Decoupling of the Homogeneous and Inhomogeneous Parts
in Inhomogeneous Quantum Groups}, 
Preprint 00-31 DMA, Universit\`a di Napoli, 
DSF/3-2001 and math.QA/0101218.

\bibitem{FioSteWes00}
G. Fiore, H. Steinacker and J. Wess
``Unbraiding the braided tensor product'',  math/0007174.

\bibitem{fuzzyq} H. Grosse, J. Madore, H. Steinacker,
  ``Field Theory on the $q$--deformed Fuzzy Sphere'', hep-th/0005273.

\end{thebibliography}
\end{document}